\documentclass[reqno,11pt]{amsart}
\usepackage{amsfonts}

\usepackage{amsmath}
\usepackage{times}
\usepackage[tcidvi]{graphicx}
\usepackage{color}

\newtheorem{theorem}{Theorem}[section]

\newtheorem{corollary}[theorem]{Corollary}

\newtheorem{lemma}[theorem]{Lemma}

\begin{document}
\title{Continuous-time quantum walk on integer lattices and homogeneous trees%
}
\author{Vladislav Kargin}
\thanks{Department of Mathematics, Stanford University, CA 94305;
kargin@stanford.edu}
\date{May 2010}
\maketitle

\begin{center}
\textbf{Abstract}
\end{center}

\begin{quotation}
This paper is concerned with the continuous-time quantum walk on $\mathbb{Z}%
, $ $\mathbb{Z}^{d},$ and infinite homogeneous trees. By using the
generating function method, we compute the limit of the average probability
distribution for the general isotropic walk on $\mathbb{Z}$, and for
nearest-neighbor walks on $\mathbb{Z}^{d}$ and infinite homogeneous trees.
In addition, we compute the asymptotic approximation for the probability of
the return to zero at time $t$ in all these cases.
\end{quotation}

\section{Introduction}

The concept of quantum walk has its origin in the field of quantum
computation where the notion of classical random walk has been adapted to
the quantum-mechanical setting in an attempt to improve the performance of
random walk algorithms. Since its origination in the middle of $90$s, this
new concept has drawn a lot of attention in physical and mathematical
literature.

The early papers that formulated the main ideas of quantum walk are %
\cite{aharonov_davidovich_zagury93} and \cite{meyer96}. From the numerous
later papers, we would like to mention \cite{farhi_gutmann98} where the
continuous-time quantum walk was defined and %
\cite{aharonov_ambainis_kempe_vazirani01} which defined and studied the
discrete-time quantum walk on finite graphs. An introductory review of
quantum walks can be found in \cite{kempe03}. For recent developments the
reader can also consult \cite{konno08}.

In general, a quantum walk is described by a triple $\left( G,\psi
,U_{t}\right) $, where $G$ is a graph, $\psi $ is a unit-length complex
vector function on this graph, i.e., $\psi \in \mathcal{H}=L^{2}\left(
G\right) \otimes \mathbb{C}^{N}$, $\Vert \psi \Vert =1$, and $U_{t}$ is a
family of unitary operators on $\mathcal{H}$.

The interpretation is that the state of a particle at time $t$ is completely
described by function $U_{t}\psi .$ Upon measurement at time $t$, the
particle is found at vertex $v$ in state $s\in \left\{ 1,2,\ldots ,N\right\} 
$ with probability $p(v,s,t)=\left| \left( U_{t}\psi \right) \left(
v,s\right) \right| ^{2}.$ (By probability we mean a collection of
non-negative numbers $p(v,s,t)$ that sum to 1, $\sum_{v\in
V}\sum_{s=1}^{N}p(v,s,t)=1.$)

There are two types of quantum walk on graph $G$. The first type is the
discrete-time walk. Time is discrete, $t\in \mathbb{Z},$ and a step of the
quantum walk is given by a unitary transformation $U,$ so that $\psi
_{t+1}=U\psi _{t}.$ This one-step transformation $U$ has some special
properties, one of which is locality: $U_{iv,ju}=0$ if the graph-theoretical
distance between vertices $\nu $ and $u$ is sufficiently large. The
discrete-time quantum walks on $\mathbb{Z}$ and $\mathbb{Z}^{d}$ have been
studied in \cite{konno02} and \cite{grimmett_janson_scudo04} who found
that their asymptotic behavior is significantly different from the behavior
of classical random walks.

In this paper we are going to investigate continuous-time quantum walks (%
\cite{farhi_gutmann98}) on infinite graphs. We assume that function $\psi $
depends only on the position of the particle and time (hence $N=1$)$,$ that $%
t\in \mathbb{R}$, and that the evolution operators $U_{t}$ are given by the
following expression:%
\begin{equation*}
U_{t}=\exp (-iXt),
\end{equation*}%
where $X$ is a self-adjoint operator (``Hamiltonian'') on $L^{2}\left(
G\right) $ that respects the structure of\ the graph. One example of such an
operator is the discrete Laplacian of the graph. We also assume that the
initial function $\psi $ is concentrated on one of the vertices, the origin
of the walk. As usual, the probability to find a particle at vertex $v$ if
the system is measured at time $t$ is given by $\left| \psi \left(
v,t\right) \right| ^{2}.$

The continuous-time walk is not local and the relation between
continuous-time and discrete-time walks is not yet clear (\cite{childs09a}%
). However, the continuous-time has an advantage over the discrete-time
model in being more tractable analytically.

For our study, we choose the simplest infinite graphs: the integer lattices $%
\mathbb{Z}$ and $\mathbb{Z}^{d},$ and the homogeneous infinite tree $\mathbb{%
T}_{m},$ in which every vertex has valency $m.$ The continuous-time walk on $%
\mathbb{Z}$ has been previously investigated by several researchers. In
particular, by using the Fourier inversion method, Gottlieb in %
\cite{gottlieb05} derived a formula for the limit probability distribution
for a very general class of quantum walks on $\mathbb{Z}$. In our paper, an
equivalent formula for the limit is obtained by using the generating series
method instead of the Fourier inversion. The advantage of this method is
that it is easier to generalize it to the case of homogeneous trees which is
not considered in \cite{gottlieb05}. In addition, this method allows us to show that the probability distribution converges to zero faster than any polynomial in time in the regions sufficiently far away from the origin. This provides a simple large deviation estimate for quantum walks on $\mathbb{Z}$.

Let us start with graph $G=\mathbb{Z}$, let the initial function $\psi
_{l}\left( 0\right) =\delta _{0},$ where $\delta _{0}$ is the delta-function
concentrated on $0,$ and let $U_{t}=\exp \left( iXt\right) ,$ where $X$ is
an Hermitian operator on $L^{2}\left( \mathbb{Z}\right) $ We are interested
in computing 
\begin{equation*}
\psi _{l}\left( t\right) =\left\langle \delta _{l}|U_{t}\delta
_{0}\right\rangle
\end{equation*}%
We say that the quantum walk has \emph{finite support} if there exists a
constant $L,$ such that $X_{ij}=0$ if $\left| i-j\right| >L$. We call the
walk \emph{isotropic }if $X_{ij}=a_{i-j}$ and $a_{-l}=\overline{a_{l}}$,
where $a_{l}$ are certain constants. If $X_{ij}=\delta _{\left| i-j\right|
-1},$ we say that the quantum walk is \emph{nearest-neighbor}.

For a general isotropic quantum walk on $\mathbb{Z}$, define the \emph{%
generating function} of the walk by the formula $P\left( \theta \right)
=\sum_{-L}^{L}a_{l}e^{il\theta }$. For example, for the nearest-neighbor
random walk, $P\left( \theta \right) =2\cos \theta .$

\begin{theorem}
\label{theorem_asymptotics_integers} Let $P\left( \theta \right) $ be the
generating function of an isotropic quantum walk with finite support on $%
\mathbb{Z}$. Assume that $\psi _{l}\left( 0\right) =\delta _{0}.$ Let $%
\alpha =l/t$ and suppose that equation $P^{\prime }\left( \theta \right)
=-\alpha $ has $K>0$ real solutions $\theta _{k}\left( \alpha \right) $ in
the interval $\left[ 0,2\pi \right) .$ Then, the transition amplitude from $%
0 $ to $l$ is given by the formula 
\begin{equation}
\psi _{l}\left( t\right) =-\frac{1}{\sqrt{2\pi t}}\sum_{k=1}^{K}\frac{1}{%
\sqrt{\left| P^{\prime \prime }\left( \theta _{k}\right) \right| }}%
e^{it\left( P\left( \theta _{k}\right) +\alpha \theta _{k}\right) \pm \pi
i/4}+O\left( \frac{1}{t}\right) ,
\end{equation}%
where the sign before $\pi i/4$ equals the sign of $P^{\prime \prime }\left(
\theta _{k}\right) .$ If equation $P^{\prime }\left( \theta \right) =-\alpha 
$ has no real solutions, then 
\begin{equation*}
\left| \psi _{l}\left( t\right) \right| \leq c_{n}t^{-n}
\end{equation*}%
for all $n$ and appropriately chosen $c_{n}.$
\end{theorem}

We will prove this theorem in Section \ref{section_QW_integers}.

Next, define the rescaled probability distribution as $p\left( \alpha
,t\right) =t\left| \psi _{\left[ \alpha t\right] }\left( t\right) \right|
^{2}$, and define the average rescaled probability as 
\begin{equation}
\overline{p}\left( \alpha ,T\right) =\frac{1}{\sqrt{T}}\int_{T}^{T+\sqrt{T}%
}p\left( \alpha ,t\right) dt  \label{formula_average_probability_definition1}
\end{equation}

\begin{corollary}
\label{corollary_probability_integers}Suppose that equation $P^{\prime
}\left( \theta \right) =-\alpha $ has $K>0$ real solutions $\theta
_{k}\left( \alpha \right) $ in the interval $\left[ 0,2\pi \right) ,$ and
that numbers 
\begin{equation*}
\omega _{k}:=P\left( \theta _{k}\right) +\alpha \theta _{k}
\end{equation*}%
are all different. Then 
\begin{equation}
\lim_{T\rightarrow \infty }\overline{p}\left( \alpha ,T\right) =\frac{1}{%
2\pi }\sum\limits_{k=1}^{K}\frac{1}{\left| P^{\prime \prime }\left( \theta
_{k}\right) \right| }.  \label{limit_distribution}
\end{equation}%
If equation $P^{\prime }\left( \theta \right) =-\alpha $ has no real
solutions, then 
\begin{equation*}
\lim_{T\rightarrow \infty }\overline{p}\left( \alpha ,T\right) =0.
\end{equation*}
\end{corollary}

In the case, when some of $\omega _{k}$ coincide, formula (\ref%
{limit_distribution}) needs a small adjustment which takes into account the
positive interference of the exponents with the same frequency.

\textbf{Proof of Corollary \ref{corollary_probability_integers}:} From
Theorem \ref{theorem_asymptotics_integers}, it follows that%
\begin{eqnarray*}
p\left( \alpha ,t\right) &=&\frac{1}{2\pi }\sum_{k=1}^{K}\frac{1}{\left|
P^{\prime \prime }\left( \theta _{k}\right) \right| } \\
&&+\frac{1}{2\pi }\sum_{k\neq l}\frac{\exp \left[ it\left( \omega
_{k}-\omega _{l}\right) +\left[ \mathrm{sgn}\left( P^{\prime \prime }\left(
\theta _{k}\right) \right) -\mathrm{sgn}\left( P^{\prime \prime }\left(
\theta _{l}\right) \right) \right] \pi i/4\right] }{\sqrt{\left| P^{\prime
\prime }\left( \theta _{k}\right) \right| \left| P^{\prime \prime }\left(
\theta _{l}\right) \right| }} \\
&&+O\left( \frac{1}{t}\right) .
\end{eqnarray*}

After averaging over $t,$ we get 
\begin{equation*}
\overline{p}\left( \alpha ,T\right) =\frac{1}{2\pi }\sum_{k=1}^{K}\frac{1}{%
\left| P^{\prime \prime }\left( \theta _{k}\right) \right| }+O\left( \frac{1%
}{\sqrt{T}}\right) .
\end{equation*}%
QED.

In particular, we obtain the following formula for the probability of return
to zero.

\begin{corollary}
Let equation $P^{\prime }\left( \theta \right) =0$ have $K>0$ real solutions 
$\theta _{k}$ in the interval $\left[ 0,2\pi \right) $ and assume that $%
P\left( \theta _{k}\right) $ are all different. Then the average probability
of the return to zero is 
\begin{equation*}
\overline{p_{0}}\left( t\right) =\frac{1}{2\pi t}\sum_{k=1}^{K}\frac{1}{%
\left| P^{\prime \prime }\left( \theta _{k}\right) \right| }+O\left( \frac{1%
}{t^{3/2}}\right)
\end{equation*}
\end{corollary}

A somewhat different expression for the limit of the average rescaled
probability distribution was obtained in \cite{gottlieb05}. Namely, by
using the Fourier inversion methods Gottlieb obtained the following formula
for the limit:%
\begin{equation*}
p_{\psi }\left( X\right) =\frac{1}{2\pi }\int_{\left\{ \theta :-P^{\prime
}\left( \theta \right) \in X\right\} }\left| \left( \mathcal{F}^{\ast }\psi
\right) \left( \theta \right) \right| ^{2}d\theta .
\end{equation*}%
Here, $p_{\psi }\left( X\right) $ denotes the limit of the probability to
find the particle in the set $\left\{ l;l/t\in X\right\} $ at time $t,$
provided that the initial wave function is $\psi $. The function $\left( 
\mathcal{F}^{\ast }\psi \right) \left( \theta \right) $ is defined as $%
\sum_{n\in Z}\psi _{n}e^{in\theta }.$ In our case, the random walk is
started from the origin and therefore $\left( \mathcal{F}^{\ast }\psi
\right) \left( \theta \right) =1.$ Then, by using Gottlieb's formula, it is
easy to compute 
\begin{equation*}
\lim_{d\alpha }\frac{p_{\psi }\left( \left[ \alpha ,\alpha +d\alpha \right]
\right) }{d\alpha }=\frac{1}{2\pi }\sum\limits_{\left\{ \theta :-P^{\prime
}\left( \theta \right) =\alpha \right\} }\frac{1}{\left| P^{\prime \prime
}\left( \theta _{k}\right) \right| },
\end{equation*}%
which is our formula (\ref{limit_distribution}). We obtain formula (\ref%
{limit_distribution}) by the method of generating functions which we will also use for quantum walks on trees.

For a particular example, consider the nearest-neighbor quantum walk on $%
\mathbb{Z}$. That is, assume that $X_{ij}=\delta _{\left| i-j\right| -1}.$
Then, $P^{\prime }\left( \theta \right) =-2\sin \theta $ and there are two
solutions $\theta _{k}\left( \alpha \right) $ for $\alpha <2$: $\theta
_{1}=\arcsin \left( \alpha /2\right) $ and $\theta _{2}=\pi -\theta _{1}.$
Hence, we can compute $\left| P^{\prime \prime }\left( \theta _{k}\right)
\right| =\sqrt{4-\alpha ^{2}},$ and the average rescaled distribution has
the following limit: 
\begin{equation}
\overline{p}\left( \alpha ,T\right) \underset{T\rightarrow \infty }{%
\rightarrow }\frac{1}{\pi }\frac{1}{\sqrt{4-\alpha ^{2}}}.
\label{arcsine_distribution}
\end{equation}%
In other words, the average rescaled distribution converges to the arcsine
law. (This result is the first limit theorem for continuous-time quantum
walks, which was obtained in \cite{konno05} by using the asymptotics of
Bessel functions).

\begin{figure}[tbph]
\begin{center}
\includegraphics[width=8cm]{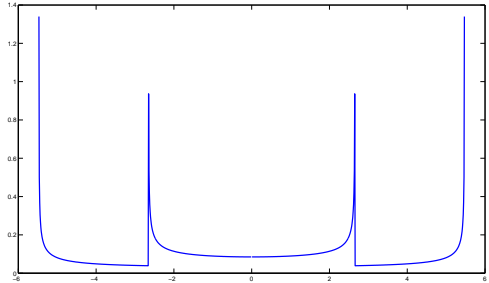}
\end{center}
\caption{The limit average probability distribution for the quantum walk on $%
\mathbb{Z}$ with $P\left( \protect\theta \right) =e^{-i2\protect\theta %
}+e^{-i\protect\theta }+e^{i\protect\theta }+e^{i2\protect\theta }.$}
\label{figure_QW_distribution_integers}
\end{figure}

Consider another example with $P\left( \theta \right) =e^{-i2\theta
}+e^{-i\theta }+e^{i\theta }+e^{i2\theta }.$ The limit probability
distribution can be computed numerically by using Theorem \ref%
{theorem_asymptotics_integers}. It is shown in Figure \ref%
{figure_QW_distribution_integers}. Note that the points where the
probability distribution has singularities correspond to local maxima of the
function $P^{\prime }\left( \theta \right) =-2(2\sin 2\theta +\sin \theta ).$

Let us summarize the main features of quantum walks on the integer lattice.
First of all, the length of the interval where the distribution of quantum
walk is essentially supported is of order $t$ instead of $\sqrt{t}$ as in
the classical case. Then, the probability of the return to the origin at
time $t$ is of order $t^{-1}$ instead of $t^{-1/2}.$ Finally, unlike the
classical case, the limit distribution is not Gaussian and its shape depends
on the generator of the quantum walk.

What can be said about the continuous-time quantum walk on $\mathbb{Z}^{d}$?
We consider here only the nearest-neighbor walk. It turns out that in this
case the quantum walk on $\mathbb{Z}^{d}$ factorizes provided it was started
from the origin. That is, every transition amplitude of the walk on $\mathbb{%
Z}^{d}$ can be written as a product of transition amplitudes of the walk on $%
\mathbb{Z}$.

For simplicity of notation we consider only the case of $\mathbb{Z}^{2}.$
The general case is similar. Let $\psi _{\left( i,j\right) }\left( t\right) $
be the transition amplitude of the transition from vertex $\left( 0,0\right) 
$ to vertex $\left( i,j\right) .$ (That is, $\psi _{\left( i,j\right)
}\left( t\right) =\left\langle \delta _{\left( i,j\right) }|U_{t}\delta
_{\left( 0,0\right) }\right\rangle ,$ where $\delta _{(0,0)}$ and $\delta
_{\left( i,j\right) }$ denote the delta-functions concentrated on vertices $%
(0,0)$ and $(i,j),$ respectively.)

\begin{theorem}
\label{theorem_factorization} Let the nearest-neighbor quantum walk on $%
\mathbb{Z}^{2}$ be started from the origin, i.e., $\psi \left( 0\right)
=\delta _{\left( 0,0\right) }$. Then, 
\begin{equation*}
\psi _{\left( i,j\right) }\left( t\right) =\psi _{i}\left( t\right) \psi
_{j}\left( t\right) ,
\end{equation*}%
where $\psi _{i}\left( t\right) $ is the transition amplitude for the
nearest-neighbor quantum walk on $\mathbb{Z}$ started from the origin.
\end{theorem}

We prove this result in \ref{section_QW_integers}. Previously, this fact was
observed without proof in Appendix of \cite{agliari_blumen_mulken08}.

Let us now turn to continuous-time quantum walks on homogeneous trees. (The
previous studies of this topic include \cite{farhi_gutmann98} and %
\cite{mulken_bierbaum_blumen06}.) We restrict our investigations to the
case of the nearest-neighbor walk.

Let $G=\mathbb{T}_{m},$ the $m$-valent infinite tree with $m\geq 3$, the
initial $\psi $ be $\delta _{e},$ where $e$ is the root of the tree, and let 
$U_{t}=\exp \left( iXt\right) ,$ where $X$ is the adjacency matrix of the
tree. Let $r:=2\sqrt{m-1}.$ (This is the spectral radius of the operator $X.$%
) Finally, let us define the following functions of parameter $\alpha $: 
\begin{eqnarray*}
\omega _{1}\left( \alpha \right) &=&\alpha \arctan \frac{\alpha }{\sqrt{%
r^{2}-\alpha ^{2}}}+\sqrt{r^{2}-\alpha ^{2}}, \\
\omega _{2}\left( \alpha \right) &=&\alpha \pi -\omega _{1},
\end{eqnarray*}%
and 
\begin{eqnarray*}
\varphi _{1}\left( \alpha \right) &=&-\arctan \left[ \frac{m}{m-2}\frac{%
\alpha }{\sqrt{r^{2}-\alpha ^{2}}}\right] -\frac{\pi }{4}, \\
\varphi _{2}\left( \alpha \right) &=&-\pi -\varphi _{1}.
\end{eqnarray*}

\begin{theorem}
\label{theorem_asymptotics_trees}Consider the nearest-neighbor quantum walk
on the regular infinite tree of valency $m.$ Assume $\psi \left( 0\right)
=\delta _{e}$, and let $\psi _{l}\left( t\right) $ be the amplitude of
transition from the root to a vertex $w$ which is located at distance $l$
from the root. Let $\alpha =l/t.$ Then%
\begin{eqnarray*}
e^{\frac{\alpha t}{2}\log \left( m-1\right) }\psi _{l}\left( t\right) &=&%
\frac{1}{\sqrt{2\pi t}}\frac{1}{\left( r^{2}-\alpha ^{2}\right) ^{1/4}}\sqrt{%
\frac{\left( m-1\right) \alpha ^{2}}{\alpha ^{2}+\left( m-2\right) ^{2}}}%
\left[ \sum_{k=1}^{2}e^{it\omega _{k}\left( \alpha \right) +i\varphi
_{k}\left( \alpha \right) }\right] \\
&&+O\left( \frac{1}{t}\right) ,\text{ if }0<\alpha <r.
\end{eqnarray*}%
In addition, there is a constant $c>0,$ which depends only on $m,$ such that%
\begin{equation*}
e^{\frac{\alpha t}{2}\log \left( m-1\right) }\left| \psi _{l}\left( t\right)
\right| \leq \frac{c}{t}\text{ if }\alpha >r.
\end{equation*}
\end{theorem}

We will prove this theorem in Section \ref{section_QW_trees}. 

Next, define 
\begin{equation}
p(\alpha ,t)=m(m-1)^{[\alpha t]-1}\cdot t|\psi _{\lbrack \alpha t]}(t)|^{2},
\label{probability_distance_trees}
\end{equation}%
and 
\begin{equation}
\overline{p}\left( \alpha ,T\right) =\frac{1}{\sqrt{T}}\int_{T}^{T+\sqrt{T}%
}p\left( \alpha ,t\right) dt
\end{equation}%
The factor $m\left( m-1\right) ^{\left[ \alpha t\right] -1}$ in (\ref%
{probability_distance_trees}) equals the number of vertices in the tree at
the distance $\left[ \alpha t\right] \geq 1$ from the root. Intuitively, $%
\overline{p}\left( \alpha ,T\right) $ is the average probability density of
the event that we find a particle at the distance approximately $\alpha t$
from the root if we measure its position at time approximately equal to $T.$
Then, we have the following corollary of Theorem \ref%
{theorem_asymptotics_trees}.

\begin{corollary}
Suppose that $\omega _{1}\left( \alpha \right) \neq \omega _{2}\left( \alpha
\right) .$ Then$\ $%
\begin{eqnarray*}
\lim_{T\rightarrow \infty }\overline{p}\left( \alpha ,T\right) &=&\frac{1}{%
\pi }\frac{1}{\sqrt{r^{2}-\alpha ^{2}}}\frac{m\alpha ^{2}}{\alpha
^{2}+\left( m-2\right) ^{2}},\text{ if }0<\alpha <r, \\
&=&0,\text{ if }\alpha >r.
\end{eqnarray*}
\end{corollary}

The proof of this corollary is similar to the proof of Corollary \ref%
{corollary_probability_integers} and is omitted.

\begin{figure}[tbph]
\begin{center}
\includegraphics[width=7cm]{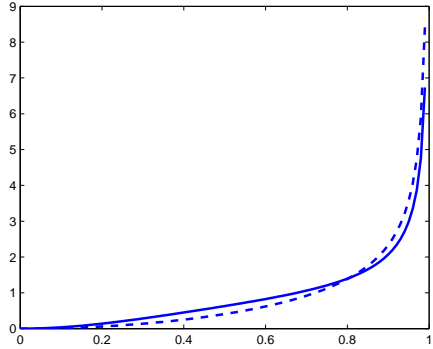}
\end{center}
\caption{The limit average probability distribution for the quantum walk on
a homogeneous tree with valency $m$. The solid line is for $m=4,$ and the
dashed line is for $m=20.$ The support of the distribution is rescaled to $%
\left[ 0,1\right] $ interval. }
\label{figure_QW_distribution_trees}
\end{figure}

The limit distribution is similar to the arcsine distribution (\ref%
{arcsine_distribution}), except it has a weighting factor, which is
different for every $m.$ A plot of the limit average distribution is shown
in Figure \ref{figure_QW_distribution_trees} for $m=4$ and $m=20.$ For the
purposes of comparison we have additionally rescaled the support of the
distribution so that supports are the same for both $m.$ It can be seen that
the walk on the tree of higher valency has higher density next to the border
of the support.

Note that the quantum walk on trees behaves (perhaps non-surprisingly) quite
different from the classical random walk on trees. In the latter case, it is
possible to show that the distribution of the particle distance from the
root is asymptotically Gaussian with mean $ct$ and standard deviation $%
\sigma \sqrt{t}$ where $c,\sigma >0$ (see \cite{sawyer_steger87} and %
\cite{lalley93}). This result is very different from what we find in the
quantum case.

\begin{figure}[tbph]
\begin{center}
\includegraphics[width=7cm]{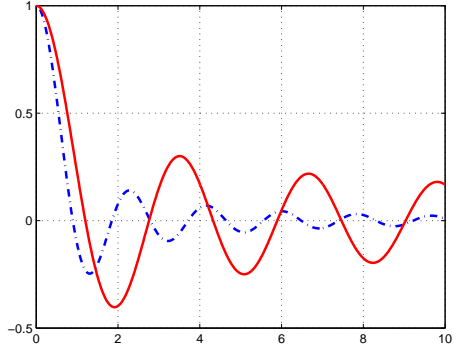}
\end{center}
\caption{The transition amplitude of the return to the root for the
nearest-neighbor quantum walk. The solid line is for integer lattice $%
\mathbb{Z}$; the dash-dotted line is for $4$-valent infinite tree $\mathbb{T}%
_{4}$.}
\label{figure_QW_return_amplitude}
\end{figure}

Consider now the amplitude of the return to zero at time $t.$

\begin{theorem}
\label{theorem_prob_return_trees}Let $\psi _{0}\left( t\right) $ denote the
transition amplitude of the return to the root at time $t$ for the
nearest-neighbor quantum walk on the $m$-valent infinite tree started at the
root$.$ Suppose that $m\geq 3.$ Then, for large $t$ the following asymptotic
approximation is valid: 
\begin{equation*}
\psi _{0}\left( t\right) =-\frac{1}{\sqrt{\pi }t^{3/2}}\frac{m\left(
m-1\right) ^{1/4}}{\left( m-2\right) ^{2}}\sin \left[ \left( 2\sqrt{m-1}%
\right) t-\pi /4\right] +O\left( \frac{1}{t^{2}}\right) .
\end{equation*}
\end{theorem}

The plot of $\psi _{0}\left( t\right) $ for $m=4$ is shown in Figure \ref%
{figure_QW_return_amplitude} by dash-dotted line. We can see that the
frequency of oscillations in the return amplitude is higher than in the case 
$m=2.$ In addition, the absolute values of maxima decline faster.

As a corollary of Theorem \ref{theorem_prob_return_trees}, we can see that
the probability of the return to zero has the following asymptotic
approximation:%
\begin{equation*}
p_{0}\left( t\right) =\frac{m^{2}\left( m-1\right) ^{1/2}}{\pi \left(
m-2\right) ^{4}}\frac{1}{t^{3}}\sin ^{2}\left[ \left( 2\sqrt{m-1}\right)
t-\pi /4\right] +O\left( t^{-7/2}\right) .
\end{equation*}

If we compare these results with the case of the classical random walk, we
find two surprising facts. First, there is no exponential decay factor in
the probability of return. The decay is polynomial of order $t^{-3}.$ Second
the exponent in this polynomial decay does not depend on the valency of the
tree although the frequency of oscillations and the overall constant does
depend on it.

The rest of the paper is organized as follows. Section \ref%
{section_QW_integers} provides proofs for Theorems \ref%
{theorem_asymptotics_integers} and \ref{theorem_factorization} concerning
quantum walks on $\mathbb{Z}$ and $\mathbb{Z}^{d}$. And Section \ref%
{section_QW_trees} gives proofs for Theorems \ref{theorem_asymptotics_trees}
and \ref{theorem_prob_return_trees} concerning the nearest-neighbour quantum
walk on homogeneous trees.

\section{Quantum walk on integer lattices}

\label{section_QW_integers}

\textbf{Proof of Theorem \ref{theorem_asymptotics_integers}:} It is
convenient to introduce additional notation. Let $G$ be a graph $\left(
V,E\right) .$ For $f,g\in L^{2}\left( V\right) ,$ we define 
\begin{equation*}
\left\langle f,g\right\rangle =\sum_{x\in V}\overline{f\left( x\right) }%
g\left( x\right) .
\end{equation*}%
Then, for $V=\mathbb{Z}$ we can write 
\begin{equation*}
\psi _{l}\left( t\right) =\left\langle \delta _{l},e^{itX}\delta
_{0}\right\rangle =\sum_{k=0}^{\infty }\left( X^{k}\right) _{0l}\frac{\left(
it\right) ^{k}}{k!}
\end{equation*}%
where $S_{-l}$ is the shift operator that sends $\delta _{k}$ to $\delta
_{k-l}.$

Let $a_{l}^{\left( k\right) }:=\left( X^{k}\right) _{0l}.$ Then, it is easy
to see that $\sum_{l=-\infty }^{\infty }a_{l}^{\left( k\right) }z^{l}=\phi
\left( z\right) ^{k},$with $\phi \left( z\right) :=\sum_{-L}^{L}a_{l}z^{l}$.
Therefore, 
\begin{equation*}
a_{l}^{\left( k\right) }=\frac{1}{2\pi i}\oint \frac{\phi \left( z\right)
^{k}}{z^{l+1}}dz,
\end{equation*}%
where the integration is over a small contour around $0.$ Hence, the
transition amplitude is given by the formula 
\begin{eqnarray*}
\psi _{l}\left( t\right) &=&\sum_{k=0}^{\infty }a_{l}^{\left( k\right) }%
\frac{\left( it\right) ^{k}}{k!} \\
&=&\frac{1}{2\pi i}\oint \frac{e^{it\phi \left( z\right) }}{z^{l+1}}dz \\
&=&-\frac{1}{2\pi }\int_{0}^{2\pi }e^{it\left( P\left( \theta \right)
+\alpha \theta \right) }d\theta ,
\end{eqnarray*}%
where we made the change of variables $z=e^{-i\theta }$ and where $\alpha
=l/t.$

We can evaluate the asymptotic behavior of this integral by the method of
stationary phase (Chapter 4 in \cite{copson67}). The points of the
stationary phase can be found from the equation: 
\begin{equation}
P^{\prime }\left( \theta \right) =-\alpha .
\label{equation_stationary_points}
\end{equation}%
Suppose that this equations has $K>0$ solutions $\theta _{k}\left( \alpha
\right) .$ Then, the asymptotic contribution of the stationary point $\theta
_{k}$ to the integral above is given by the following formula: 
\begin{equation*}
-\frac{1}{\sqrt{2\pi t}}\frac{1}{\sqrt{\left| P^{\prime \prime }\left(
\theta _{k}\right) \right| }}e^{it\left( P\left( \theta _{k}\right) +\alpha
\theta _{k}\right) \pm \pi i/4}+O\left( \frac{1}{t}\right) ,
\end{equation*}%
where the sign before $\pi i/4$ depends on whether $P^{\prime \prime }\left(
\theta _{k}\right) $ is positive or negative. By adding these contributions
we obtain the first claim of the theorem.

If the equation (\ref{equation_stationary_points}) has no real solutions,
then there are no points of stationary phase in interval $\left[ 0,2\pi
\right) .$ In this case, we can apply the method of integration by parts.
Usually, in this case the asymptotic approximation is of the order $t^{-1}.$
However, in our case it is smaller due to special properties of function $%
P\left( \theta \right) $ and number $\alpha .$

Indeed, let us denote $P\left( \theta \right) +\alpha \theta $ as $f_{\alpha
}\left( \theta \right) $ for shortness. Note that the first derivative $%
f_{\alpha }^{\prime }$ is periodic with period $2\pi $ and therefore all
other derivatives are also periodic with period $2\pi .$ In addition, if $%
\alpha t$ is integer (which is exactly the case we consider, then $\exp
\left( itf_{\alpha }\left( \theta \right) \right) $ is periodic with period $%
2\pi .$

Since $\left| f^{\prime }\left( \theta \right) \right| \neq 0$ anywhere in
interval $\left[ 0,2\pi \right) ,$ hence we can use integration by parts in
the following form: 
\begin{eqnarray*}
\int_{0}^{2\pi }e^{itf_{\alpha }\left( \theta \right) }d\theta &=&\frac{1}{it%
}\int_{0}^{2\pi }\frac{1}{f^{\prime }\left( \theta \right) }\frac{d}{d\theta 
}\left( e^{itf_{\alpha }\left( \theta \right) }\right) d\theta \\
&=&\frac{1}{it}\left[ \frac{1}{f^{\prime }\left( 2\pi \right) }\left(
e^{itf_{\alpha }\left( 2\pi \right) }\right) -\frac{1}{f^{\prime }\left(
0\right) }\left( e^{itf_{\alpha }\left( 0\right) }\right) \right] \\
&&-\frac{1}{it}\int_{0}^{2\pi }\frac{d}{d\theta }\left( \frac{1}{f^{\prime
}\left( \theta \right) }\right) e^{itf_{\alpha }\left( \theta \right)
}d\theta .
\end{eqnarray*}%
By using the special properties of the function $f_{\alpha }\left( \theta
\right) ,$ we can conclude that the first term is zero and therefore 
\begin{equation*}
\int_{0}^{2\pi }e^{itf_{\alpha }\left( \theta \right) }d\theta =-\frac{1}{it}%
\int_{0}^{2\pi }\left( \frac{1}{f^{\prime }\left( \theta \right) }\right)
^{\prime }e^{itf_{\alpha }\left( \theta \right) }d\theta .
\end{equation*}%
In particular, this integral is $O\left( t^{-1}\right) .$ Since the function 
$\left( 1/f^{\prime }\left( \theta \right) \right) ^{\prime }$ is periodic,
hence the argument can be repeated . It is easy to see that it can be
repeated indefinitely, and we obtain that the integral is less than $%
c_{n}t^{-n}$ for every $n.$ QED.

\textbf{Proof of Theorem \ref{theorem_factorization}:} Let $X$ be the
adjacency matrix for $\mathbb{Z}^{2}$ and $H$ and $V$ be the adjacency
matrices that take into account only horizontal and vertical bonds,
respectively. In other words, 
\begin{equation*}
X=H+V,
\end{equation*}

and 
\begin{eqnarray*}
H_{ij,kl} &=&\left( \delta _{i-1,k}+\delta _{i+1,k}\right) \delta _{jl}, \\
V_{ij,kl} &=&\delta _{ik}\left( \delta _{j-1,l}+\delta _{j+1,l}\right) .
\end{eqnarray*}%
It is easy to see that 
\begin{equation*}
\left( HV\right) _{ij,kl}=\left( VH\right) _{ij,kl}=\left( \delta
_{i-1,k}+\delta _{i+1,k}\right) \left( \delta _{j-1,l}+\delta
_{j+1,l}\right) .
\end{equation*}%
That is, $H$ and $V$ commute. This implies that 
\begin{equation*}
e^{itX}=e^{itH}e^{itV}.
\end{equation*}

After we apply $e^{itV}$ to $\psi \left( 0\right) =\delta _{\left(
0,0\right) },$ the result at vertex $\left( i,j\right) $ is $\delta
_{0i}\psi _{j}\left( t\right) $ where $\delta _{0i}$ is the Kronecker delta
and $\psi _{j}\left( t\right) $ is the wave function for the
nearest-neighbor quantum walk at $\mathbb{Z}$ started at $0.$ Next, after we
apply $e^{itH},$ the result at vertex $i,j$ is $\psi _{i}\left( t\right)
\psi _{j}\left( t\right) $ since it equals the wave function of the
nearest-neighbor walk on graph $\mathbb{Z\times }\left( 0,j\right) $ started
with the initial data $\psi _{j}\left( t\right) .$ QED.

\section{Nearest-neighbor quantum walk on trees}

\label{section_QW_trees}

\textbf{Proof of Theorem \ref{theorem_asymptotics_trees}:}

First, note that 
\begin{eqnarray*}
\psi _{w}\left( t\right) &=&\left\langle \delta _{w},e^{itX}\delta
_{e}\right\rangle =\sum_{k=0}^{\infty }\left\langle \delta _{w},X^{k}\delta
_{e}\right\rangle \frac{\left( it\right) ^{k}}{k!} \\
&&\sum_{k=0}^{\infty }c_{k}\left( \left| w\right| \right) \frac{\left(
it\right) ^{k}}{k!},
\end{eqnarray*}%
where $c_{k}\left( \left| w\right| \right) $ denotes the number of all
possible paths with $k$ edges that start at the root and end at vertex $w$.

Let $A_{k}$ denote the number of paths from $e$ to $e$ that have length $k$
and do not pass along a specific edge which is connected to $e,$ say, do not
pass along edge $x_{1}.$ Let $B_{k}$ be the number of paths from $e$ to $e$
that have length $k,$ without any additional restrictions. Let $A\left(
z\right) $ and $B\left( z\right) $ denote the generating functions for $%
A_{k} $ and $B_{k},$ respectively, that is, 
\begin{equation*}
A\left( z\right) =\sum_{k=0}^{\infty }A_{k}z^{k},\text{ and }B\left(
z\right) =\sum_{k=0}^{\infty }B_{k}z^{k},
\end{equation*}%
where we set $A_{0}=B_{0}=1.$

From Lemma \ref{lemma_ckl} proved below, it follows that 
\begin{equation*}
\sum_{r=0}^{\infty }c_{l+r}\left( l\right) z^{r}=A\left( z\right)
^{l}B\left( z\right) .
\end{equation*}%
Hence, 
\begin{equation*}
c_{l+r}\left( l\right) =\frac{1}{2\pi i}\oint \frac{A\left( z\right)
^{l}B\left( z\right) }{z^{r+1}}dz,
\end{equation*}%
where the integration is around a small circle around $0.$

Let $\left| w\right| =l,$ then for the transition amplitude from $e$ to $w$,
we can write%
\begin{eqnarray*}
\psi _{l}\left( t\right) &=&\sum_{k=l}^{\infty }\frac{\left( it\right) ^{k}}{%
k!}c_{k}\left( l\right) \\
&&\sum_{r=0}^{\infty }\frac{\left( it\right) ^{l+r}}{\left( l+r\right) !}%
\frac{1}{2\pi i}\oint \frac{A\left( z\right) ^{l}B\left( z\right) }{z^{r+1}}%
dz,
\end{eqnarray*}%
which we can re-write as follows: 
\begin{eqnarray*}
\psi _{l}\left( t\right) &=&\frac{1}{2\pi i}\oint A\left( z\right)
^{l}B\left( z\right) z^{l-1}\left( \sum_{k=0}^{\infty }\frac{\left(
it/z\right) ^{l+k}}{\left( l+k\right) !}\right) dz \\
&=&\frac{1}{2\pi i}\oint A\left( z\right) ^{l}B\left( z\right) z^{l-1}\left[
e^{it/z}-\sum_{k=0}^{l-1}\frac{\left( it/z\right) ^{k}}{k!}\right] dz.
\end{eqnarray*}%
The sum in the last line gives zero contribution to the integral since
neither $A\left( z\right) $ nor $B\left( z\right) $ has any singularity at $%
0.$ Hence, we can write%
\begin{eqnarray*}
\psi _{l}\left( t\right) &=&\frac{1}{2\pi i}\oint A\left( z\right)
^{l}B\left( z\right) z^{l-1}e^{it/z}dz \\
&=&\frac{1}{2\pi i}\oint \left[ \frac{A\left( 1/u\right) }{u}\right] ^{l}%
\frac{B\left( 1/u\right) }{u}e^{itu}du \\
&=&\frac{1}{2\pi i}\oint \left[ F\left( u\right) \right] ^{l}G(u)e^{itu}du,
\end{eqnarray*}%
where we used the substitution $u=1/z,$ and 
\begin{eqnarray*}
F\left( u\right) &:&=\frac{1}{u}A\left( \frac{1}{u}\right) , \\
G\left( u\right) &:&=\frac{1}{u}B\left( \frac{1}{u}\right) .
\end{eqnarray*}%
The second and third integrals are taken over a sufficiently large circle
around the zero which includes all of the singularities of $F\left( u\right) 
$ and $G\left( u\right) .$

We calculate $F\left( u\right) $ and $G\left( u\right) $ explicitly below
(Lemma \ref{lemma_G} and \ref{lemma_F}). The function $G\left( u\right) $ is
analytical at points $u=\pm m,$ therefore the only singularities of the
integrand are branch points of $F(u)$ and $G(u)$ at $u=\pm 2\sqrt{m-1}.$

We want to find out the asymptotic approximation for those values of $l$
which are comparable with $t.$ Let $l=\alpha t$ with $\alpha \geq 0.$ The we
can write the transition amplitude as follows:

\begin{equation}
\psi _{l}\left( t\right) =\frac{1}{2\pi i}\oint e^{it[u-i\alpha \log
F(u)]}G(u)du.  \label{transition_amplitude}
\end{equation}%
Recall that $r:=2\sqrt{m-1}.$ Let us deform the contour of integration so
that it goes first from $-r$ to $r$ just below the real axis, and then goes
back just above the real axis.

Let $f\left( u\right) :=u-i\alpha \log F(u).$ For real $u\in \left[ -r,r%
\right] ,$ we can compute $\mathrm{Im}f\left( u\right) =(\alpha /2)\log
\left( m-1\right) ,$ which is constant with respect to $u.$ Hence, we can
use the stationary phase approximation to this integral.

In order to find the points of stationary phase, we need to solve the
equation $d\left( \mathrm{Re}f\left( u\right) \right) /du=0.$ Since $\mathrm{%
Im}f\left( u\right) $ is constant, it is the same as solving 
\begin{equation}
\frac{df\left( u\right) }{du}\equiv 1-i\alpha \frac{F^{\prime }\left(
u\right) }{F\left( u\right) }=0.  \label{equation_stationary_points_trees}
\end{equation}

First, let us consider the case $0<\alpha <r$.

For the part of the contour that lies in the upper part of the complex
plane, we have: $F\left( u\right) =\left( u-i\sqrt{r^{2}-u^{2}}\right)
/\left( 2\left( m-1\right) \right) ,$ hence $f^{\prime }\left( u\right)
=1+\alpha /\sqrt{r^{2}-u^{2}}$ and equation (\ref%
{equation_stationary_points_trees}) becomes $-\alpha =\sqrt{r^{2}-u^{2}}$
which has no solutions in the interval $\left( -r,r\right) $ for any
positive $\alpha .$ Hence, the contribution of this part of the contour is
asymptotically negligible provided that the integral along the other part of
the contour has stationary points.

For the part of the contour that lies in the lower part of the complex
plane, we have $F\left( u\right) =\left( u+i\sqrt{r^{2}-u^{2}}\right)
/\left( 2\left( m-1\right) \right) ,$ and the equation (\ref%
{equation_stationary_points_trees}) reduces to $\alpha =\sqrt{r^{2}-u^{2}},$
which has two solutions $u_{1,2}=\pm \sqrt{r^{2}-\alpha ^{2}}$ for $\alpha
<r.$

Recall that the method of stationary phase says that if $\overline{u}$ is
the only stationary point of function $f\left( u\right) ,$ located inside $%
\left[ a,b\right] ,$ then 
\begin{equation*}
\int_{a}^{b}e^{itf\left( u\right) }G\left( u\right) du=\sqrt{\frac{2\pi }{%
tf^{\prime \prime }\left( \overline{u}\right) }}G\left( \overline{u}\right)
e^{itf\left( \overline{u}\right) \pm \pi i/4}+O\left( 1/t\right) ,
\end{equation*}%
where the sign before $\pi i/4$ is positive if $\mathrm{Re}f^{\prime \prime
}\left( \overline{u}\right) >0$ and negative if $\mathrm{Re}f^{\prime \prime
}\left( \overline{u}\right) <0$.

We compute $F\left( u_{1,2}\right) =\left( \pm \sqrt{r^{2}-\alpha ^{2}}%
+i\alpha \right) /\left( 2\left( m-1\right) \right) .$ The second derivative
of $f\left( u\right) $ can be evaluated at $u_{1,2}$ as follows. 
\begin{equation*}
f^{\prime \prime }\left( u_{1,2}\right) =\mp \frac{\sqrt{r^{2}-\alpha ^{2}}}{%
\alpha ^{2}}.
\end{equation*}

In addition, we have 
\begin{equation*}
G\left( u_{1,2}\right) =\frac{\pm \left( m-2\right) \sqrt{r^{2}-\alpha ^{2}}%
-im\alpha }{2\left( \alpha ^{2}+(m-2)^{2}\right) },
\end{equation*}%
and 
\begin{equation*}
\left| G\left( u_{1,2}\right) \right| =\sqrt{\frac{m-1}{\alpha ^{2}+\left(
m-2\right) ^{2}}}
\end{equation*}

Hence, \ 
\begin{eqnarray*}
\psi _{l}\left( t\right) &=&e^{-\frac{\alpha t}{2}\log \left( m-1\right)
}\{\left( \frac{\alpha ^{2}}{2\pi t}\frac{1}{\left( r^{2}-\alpha ^{2}\right)
^{1/2}}\right) ^{1/2}\sqrt{\frac{m-1}{\alpha ^{2}+\left( m-2\right) ^{2}}} \\
&&\times \left[ \sum_{k=1}^{2}e^{it\omega _{k}\left( \alpha \right)
+i\varphi _{k}\left( \alpha \right) }\right] +O\left( \frac{1}{t}\right) \}.
\end{eqnarray*}

Here, the frequencies can be computed as 
\begin{equation*}
\omega _{1}=\alpha \arctan \frac{\alpha }{\sqrt{r^{2}-\alpha ^{2}}}+\sqrt{%
r^{2}-\alpha ^{2}},
\end{equation*}%
and 
\begin{equation*}
\omega _{2}=\alpha \pi -\omega _{1},
\end{equation*}%
and the phases can be computed as 
\begin{equation*}
\varphi _{1}=-\arctan \frac{m\alpha }{\left( m-2\right) \sqrt{r^{2}-\alpha
^{2}}}-\frac{\pi }{4},
\end{equation*}%
and\ 
\begin{equation*}
\varphi _{2}=-\pi -\varphi _{1}.
\end{equation*}

Now, consider the case $\alpha >r.$ In this case, neither part of the
contour has a point of stationary phase and for the large $t$, the boundary
points of the interval $\left[ -r,r\right] $ contribute most to the
integral. In this situation, we can estimate the integral by using
integration by parts. Consider, for example, the integral 
\begin{equation*}
I_{1}\left( \varepsilon \right) =\int_{-r+\varepsilon }^{r-\varepsilon
}e^{itf\left( u\right) }G\left( u\right) du,
\end{equation*}%
where $f\left( u\right) $ and $G\left( u\right) $ are defined as continuous
limits of the upper half-plane branches of $f\left( u\right) $ and $G\left(
u\right) .$ Then we can write: 
\begin{eqnarray}
I_{1}\left( \varepsilon \right) &=&\frac{1}{it}\int_{-r+\varepsilon
}^{r-\varepsilon }\left( \frac{d}{du}e^{itf\left( u\right) }\right) \frac{%
G\left( u\right) }{f^{\prime }\left( u\right) }du  \notag \\
&=&\left. e^{itf\left( u\right) }\frac{G\left( u\right) }{itf^{\prime
}\left( u\right) }\right| _{-r+\varepsilon }^{r+\varepsilon }-\frac{1}{it}%
\int_{-r+\varepsilon }^{r-\varepsilon }e^{itf\left( u\right) }\left( \frac{%
G\left( u\right) }{f^{\prime }\left( u\right) }\right) ^{\prime }du.
\label{formula_integral_by_parts}
\end{eqnarray}

Since $f^{\prime }\left( u\right) =1+\alpha /\sqrt{r^{2}-u^{2}},$ therefore $%
\frac{G\left( u\right) }{f^{\prime }\left( u\right) }\rightarrow 0$ as $%
u\rightarrow \pm r,$ which implies that the first part of (\ref%
{formula_integral_by_parts}) becomes zero as $\varepsilon \rightarrow 0.$
For the second part, note that 
\begin{equation*}
G^{\prime }\left( u\right) =A(u)+B(u)\sqrt{r^{2}-u^{2}}+C\left( u\right) 
\frac{1}{\sqrt{r^{2}-u^{2}}},
\end{equation*}%
for some functions $A\left( u\right) ,$ $B\left( u\right) ,$ and $C\left(
u\right) $ analytic on $\left[ -r,r\right] ,$ which implies that 
\begin{equation*}
\left( \frac{G\left( u\right) }{f^{\prime }\left( u\right) }\right) ^{\prime
}=\frac{G^{\prime }\left( u\right) }{f^{\prime }\left( u\right) }-\frac{%
G\left( u\right) }{\left[ f^{\prime }\left( u\right) \right] ^{2}}%
f^{^{\prime \prime }}
\end{equation*}%
has singularities $\left( u+r\right) ^{-1/2}$ and $\left( r-u\right) ^{-1/2}$
at $-r$ and $r$ respectively$.$ This implies that $\left( \frac{G\left(
u\right) }{f^{\prime }\left( u\right) }\right) ^{\prime }$ is absolutely
integrable at $\left[ -r,r\right] $ and therefore%
\begin{eqnarray*}
\left| \lim_{\varepsilon \rightarrow 0}I_{1}\left( \varepsilon \right)
\right| &\leq &\frac{1}{t}e^{-t\mathrm{Im}f}\int_{-r}^{r}\left| \left( \frac{%
G\left( u\right) }{f^{\prime }\left( u\right) }\right) ^{\prime }\right| du
\\
&\leq &\frac{c}{t}e^{-\left( \alpha t/2\right) \log \left( m-1\right) }.
\end{eqnarray*}

A similar estimate holds for the integral along the contour in the lower
half-plane. This completes the proof of Theorem \ref%
{theorem_asymptotics_trees}.

Here are the auxiliary results that we used in the proof.

\begin{lemma}
\label{lemma_ckl} Suppose that graph $G$ is an infinite homogeneous tree
with root $e.$ Let $A_{k}$ be the number of paths in $G$ from $e$ to $e$
that have length $k$ and do not pass along a specific edge which is
connected to $e.$ Let $B_{k}$ be the number of paths from $e$ to $e$ that
have length $k,$ without further restrictions, and let $c_{k}\left( \left|
w\right| \right) $ be the number of paths from $e$ to $w$ that have length $%
k $. Then, 
\begin{equation*}
c_{k}\left( \left| w\right| \right) =\sum\limits_{k_{0}+k_{1}+\ldots
+k_{\left| w\right| }=k}A_{k_{0}}A_{k_{1}}\ldots A_{k_{\left| w\right|
-1}}B_{k_{\left| w\right| }}.
\end{equation*}
\end{lemma}

\textbf{Proof:} \ Assume that each edge in the tree is oriented and has a
label, $x$, which is chosen from the set $\left\{ 1,\ldots ,m\right\} .$ It
is assumed that that the labels of edges around each vertex are all
different. We write label $x$ if we move in the direction of the orientation
and $x^{-1}$ if we move in the opposite direction. Let $x_{l}x_{l-1}...x_{1}$
be the shortest path from $e$ to $w.$ $\ $There is a one-to-one
correspondence between the set of shortest paths and vertices so we can
write $w=x_{l}x_{l-1}...x_{1}.$ Also, let $w_{i}=x_{i}x_{i-1}...x_{1}.$ This
is one of the vertices on the shortest path from $e$ to $w$. We write the
edges in the path from right to left so that $w_{1}$ is a neighbor of the
root.

Every path from $e$ to $w$ can be considered as the shortest path from $e$
to $w$ decorated with loops which can be attached at each of the points of
the shortest path, $w_{i}.$ In order to make sure that we do not double
count the loops we forbid the loop attached at $w_{i}$ to go along the edge
that connects $w_{i}$ to $w_{i+1}.$ In this way, at every point of the path
we know in which loop we are in: We are always in the loop attached at that $%
w_{i}$ that has the largest length $\left| w_{i}\right| $ among all those
vertices $w_{i}$ that have already been visited.

Let $l=\left| w\right| .$ The number of possible different loops that can be
attached at $w_{0},$ $w_{1},$ \ldots , $w_{l-1}$ is counted by $A_{k_{0}},$ $%
A_{k_{1}},$ \ldots , $A_{k_{l-1}},$ respectively, where $k_{0},k_{1},$
\ldots , $k_{l-1}$ are the lengths of the loops. The number of different
loops that can be attached at $w=w_{l}$ is counted by $B_{k_{l}}.$ Then, the
total length of the path is $k_{0}+k_{1}+\ldots +k_{l}$ and by assumption it
must be equal to $k.$ Hence the total number of paths is 
\begin{equation*}
\sum\limits_{k_{0}+k_{1}+\ldots +k_{l}=k}A_{k_{0}}A_{k_{1}}\ldots
A_{k_{l-1}}B_{k_{l}}.
\end{equation*}%
QED.

\begin{lemma}
\label{lemma_G}%
\begin{equation}
G\left( z\right) :=\frac{1}{z}B\left( \frac{1}{z}\right) =\frac{-\left(
m-2\right) z+m\sqrt{z^{2}-4\left( m-1\right) }}{2\left( z^{2}-m^{2}\right) }.
\end{equation}
\end{lemma}

\textbf{Proof:} The function $B\left( z\right) $ is related to the Green
function of the nearest-neighbor random walk on an infinite tree, which is
well-known. (See Dynkin and Malyutov for the seminal contribution, and Lemma
1.24 on p. 9 in \cite{woess00}.) Hence, we can compute 
\begin{equation*}
B\left( z\right) =\frac{-\left( m-2\right) +m\sqrt{1-4\left( m-1\right) z^{2}%
}}{2\left( 1-m^{2}z^{2}\right) }.
\end{equation*}%
It follows that 
\begin{equation}
G\left( z\right) =\frac{-\left( m-2\right) z+m\sqrt{z^{2}-4\left( m-1\right) 
}}{2\left( z^{2}-m^{2}\right) }.
\end{equation}

QED

Note that we chose the branches of $G\left( z\right) $ in such a way that
the function is analytical outside the cut $\left[ -2\sqrt{m-1},2\sqrt{m-1}%
\right] $. In particular, this function does not have poles at $\pm m$.

More precisely, the sign before the square root is determined by the rule
that for sufficiently small $t,$ 
\begin{equation*}
G\left( it\right) \approx -i\frac{\sqrt{m-1}}{m}\in \mathbb{C}^{-}
\end{equation*}%
and 
\begin{equation*}
G\left( -it\right) \approx i\frac{\sqrt{m-1}}{m}\in \mathbb{C}^{+}.
\end{equation*}

\begin{lemma}
\label{lemma_F}%
\begin{equation*}
F\left( z\right) :=\frac{1}{z}A\left( \frac{1}{z}\right) =\frac{z-\sqrt{%
z^{2}-4\left( m-1\right) }}{2\left( m-1\right) }.
\end{equation*}
\end{lemma}

\textbf{Proof:} In order to compute $A\left( z\right) ,$ we note that the
following recursive relation holds. 
\begin{equation}
A_{2k}=\left( m-1\right) \sum_{l=0}^{k-1}A_{2l}A_{2\left( k-l-1\right) }.
\label{formula_A_recursion}
\end{equation}

Indeed, consider a path from $e$ to $e,$ that avoids the edge $x_{1}.$ There
are $m-1$ possibilities to start the path. Suppose that the path starts with 
$x_{i},$ $i\neq 1,$ so that the second point on the path is the endpoint of $%
x_{i}$ which we denote $w_{1}.$ Let $r$ \ be the first time when the path
returns to $e.$ Then $w_{r-1}=w_{1}$ and the path from $w_{1}$ to $w_{r-1}$
is one of the $A_{r-2}$ paths from $w_{1}$ to $w_{1}$ that avoid passing
through the edge labelled $x_{i}.$ The remainder of the path goes from $e$
to $e$ and it is one of the $A_{2k-r}$ paths that avoid the edge $x_{1}.$
The number $r$ must be even, greater than $0$ and less than $2k.$ Hence we
can write it as $r=2l+2,$ where $0\leq l\leq k-1.$ This implies the
recursive formula (\ref{formula_A_recursion}).

Next, we can use the recursion formula for Catalan numbers, 
\begin{equation*}
C_{k}=\sum_{l=0}^{k-1}C_{l}C_{k-l-1},
\end{equation*}%
and formula (\ref{formula_A_recursion}) in order to conclude that 
\begin{equation*}
A_{2k}=\left( m-1\right) ^{k}C_{k}.
\end{equation*}%
By using the generating function for Catalan numbers, we obtain the
following formula for $A\left( z\right) $: 
\begin{equation*}
A\left( z\right) =\frac{1-\sqrt{1-4\left( m-1\right) z^{2}}}{2\left(
m-1\right) z^{2}}.
\end{equation*}

It follows that 
\begin{equation*}
F\left( z\right) =\frac{z-\sqrt{z^{2}-4\left( m-1\right) }}{2\left(
m-1\right) }.
\end{equation*}%
QED.

The sign of the square root in the expression for $F\left( z\right) $ is
determined by the following rule: for all sufficiently small $t,$%
\begin{equation*}
F\left( it\right) \approx -i/\sqrt{m-1}\in \mathbb{C}^{-},
\end{equation*}%
and 
\begin{equation*}
F\left( -it\right) \approx i/\sqrt{m-1}\in \mathbb{C}^{+}.
\end{equation*}

\textbf{Proof of Theorem \ref{theorem_prob_return_trees}:} By (\ref%
{transition_amplitude}), we need to find asymptotics for 
\begin{equation}
\psi _{0}\left( t\right) =\frac{1}{2\pi i}\oint e^{itu}G(u)du,
\end{equation}%
where 
\begin{equation*}
G\left( u\right) =\frac{-\left( m-2\right) u+m\sqrt{u^{2}-r^{2}}}{2\left(
u^{2}-m^{2}\right) }.
\end{equation*}%
and $r=2\sqrt{m-1}.$ We can deform the contour so that it starts at $-r,$
passes just below the real axis to $r$ and then returns back to $-r$ just
above the real axis. Then, we find that 
\begin{equation*}
\psi _{0}\left( t\right) =\frac{1}{\pi }\int_{-r}^{r}e^{itu}\frac{m\sqrt{%
r^{2}-u^{2}}}{2\left( u^{2}-m^{2}\right) }du,
\end{equation*}

The main contribution is produced by singular points $\pm r.$ After
integration by parts, we obtain the following formula. 
\begin{equation*}
\psi _{0}\left( t\right) =-\frac{1}{2\pi i}\frac{m}{t}\left[ \int_{-r}^{r}%
\frac{1}{\sqrt{r^{2}-u^{2}}}\frac{-ue^{itu}}{u^{2}-m^{2}}du+\int_{-r}^{r}%
\frac{-2u\sqrt{r^{2}-u^{2}}}{\left( u^{2}-m^{2}\right) ^{2}}e^{itu}du\right]
.
\end{equation*}%
We can apply van der Corput's results (see \cite{copson67}, p. 24) to the
first integral in the brackets and obtain the following asymptotic
approximation.

\begin{equation*}
\int_{-r}^{r}\frac{1}{\sqrt{r^{2}-u^{2}}}\frac{-ue^{itu}}{u^{2}-m^{2}}du=%
\sqrt{\frac{\pi }{2rt}}\frac{r}{r^{2}-m^{2}}\left( e^{-itr+\pi
i/4}-e^{itr-\pi i/4}\right) +O\left( t^{=1}\right) .
\end{equation*}%
We can apply the integration by parts to the second integral and find that
it is $O\left( t^{-1}\right) .$ It follows that 
\begin{eqnarray*}
\psi _{0}\left( t\right) &=&\frac{1}{\sqrt{2\pi }t^{3/2}}\frac{m\sqrt{r}}{%
r^{2}-m^{2}}\sin \left( rt-\pi /4\right) +O\left( t^{-2}\right) \\
&=&-\frac{1}{\sqrt{\pi }t^{3/2}}\frac{m\left( m-1\right) ^{1/4}}{\left(
m-2\right) ^{2}}\sin \left( 2\sqrt{m-1}t-\pi /4\right) +O\left(
t^{-2}\right) .
\end{eqnarray*}%
QED.

\bibliographystyle{plain}
\bibliography{comtest}

\begin{thebibliography}{10}

\bibitem{agliari_blumen_mulken08}
E.~Agliari, A.~Blumen, and O.~Mulken.
\newblock Dynamics of continuous-time quantum walks in restricted geometries.
\newblock {\em Journal of Physics A: Mathematical and Theoretical}, 41:445301,
  2008.

\bibitem{aharonov_ambainis_kempe_vazirani01}
D.~Aharonov, A.~Ambainis, J.~Kempe, and U.~Vazirani.
\newblock Quantum walks on graphs.
\newblock In {\em Proceedings of the 33rd STOC}, pages 50--59. ACM, New York,
  2001.
\newblock arxiv:quant-ph/0012090v2 25 May 2002.

\bibitem{aharonov_davidovich_zagury93}
Y.~Aharonov, L.~Davidovich, and N.~Zagury.
\newblock Quantum random walks.
\newblock {\em Physics Review A}, 48:1687--1690, 1993.

\bibitem{childs09a}
Andrew~M. Childs.
\newblock On the relationship between continuous- and discrete-time quantum
  walk.
\newblock {\em Communications in Mathematical Physics}, 2009.
\newblock also available at http://www.arxiv.org/abs/0810.0312v2.

\bibitem{copson67}
E.~T. Copson.
\newblock {\em Asymptotic expansions}.
\newblock Cambridge Tracts in Mathematics and Mathematical Physics. Cambridge
  University Press, 1967.

\bibitem{farhi_gutmann98}
E.~Farhi and S.~Gutmann.
\newblock Quantum computation and decision trees.
\newblock {\em Physics Review A}, 58:915--928, 1998.

\bibitem{gottlieb05}
Alex~D. Gottlieb.
\newblock Convergence of continuous-time quantum walks on the line.
\newblock {\em Physical Review E}, 72:047102, 2005.

\bibitem{grimmett_janson_scudo04}
Geoffrey Grimmett, Svante Janson, and Petra~F. Scudo.
\newblock Weak limits for quantum random walks.
\newblock {\em Physical Review E}, 69:026119, 2004.

\bibitem{kempe03}
Julia Kempe.
\newblock Quantum random walks - an introductory overview.
\newblock {\em Contemporary Physics}, 44:302--327, 2003.
\newblock arxiv:quant-ph/0303081v1.

\bibitem{konno02}
Norio Konno.
\newblock Quantum random walks in one dimension.
\newblock {\em Quantum Information Processing}, 1:345--354, 2002.

\bibitem{konno05}
Norio Konno.
\newblock Limit theorem for continuous-time quantum walk on the line.
\newblock {\em Physical Review E}, 72:026113, 2005.

\bibitem{konno08}
Norio Konno.
\newblock Quantum walks.
\newblock In {\em Quantum Potential Theory}, volume 1954 of {\em Lecture Notes
  in Mathematics}, pages 309--452. Springer, Berlin, 2008.

\bibitem{lalley93}
Steven~P. Lalley.
\newblock Finite range random walk on free groups and homogeneous trees.
\newblock {\em Annals of Probability}, 21:2087--2130, 1993.

\bibitem{meyer96}
D.~Meyer.
\newblock From quantum cellular automata to quantum lattice gases.
\newblock {\em Journal of Statistical Physics}, 85:551--574, 1996.

\bibitem{mulken_bierbaum_blumen06}
O.~Mulken, V.~Bierbaum, and A.~Blumen.
\newblock Coherent exciton transport in dendrimers and continuous-time quantum
  walks.
\newblock {\em Journal of Chemical Physics}, 124:124905, 2006.

\bibitem{sawyer_steger87}
S.~Sawyer and T.~Steger.
\newblock The rate of escape for anisotropic random walks in a tree.
\newblock {\em Probability Theory and Related Fields}, 76:207--230, 1987.

\bibitem{woess00}
Wolfgang Woess.
\newblock {\em Random Walks on Infinite Graphs and Groups}.
\newblock Cambridge Tracts in Mathematics. Cambridge University Press, 2000.

\end{thebibliography}

\end{document}